\documentclass[12pt,reqno]{amsart}
\usepackage{cite}
\usepackage{amsfonts}
\usepackage{amssymb, amsmath, latexsym}
\usepackage{longtable}
\usepackage{lipsum}
\usepackage{hyperref}
\usepackage{mleftright}
\usepackage{graphicx,mathtools}
\usepackage{mathtools}
\DeclarePairedDelimiter\floor{\lfloor}{\rfloor}
\usepackage{array}

\newcommand\cplus{\mathbin{\raisebox{-\height}{$+$}}}
\newcommand\contdots{\raisebox{-\height}{$\vphantom{+}\dotsm$}}

\numberwithin{equation}{section}
\theoremstyle{plain}
\newtheorem{theorem}{Theorem}[section]

\newtheorem{lemma}[theorem]{Lemma}
\newtheorem{conjecture}[theorem]{Conjecture}

\newtheorem{remark}[theorem]{Remark}
\newcommand\T{\rule{0pt}{2.6ex}}
\newcommand\B{\rule[-1.2ex]{0pt}{0pt}}

\title[Congruences for $k$-elongated plane partition diamonds] {Congruences for $k$-elongated plane partition diamonds}
\author[N. D. Baruah]{Nayandeep Deka Baruah}
\address{Department of Mathematical Sciences, Tezpur University,  Assam 784028, India}
\email{nayan@tezu.ernet.in}
\author[H. Das]{Hirakjyoti Das}
\address{Department of Mathematical Sciences, Tezpur University,  Assam 784028, India}
\email{hdas@tezu.ernet.in}
\author[P. Talukdar]{Pranjal Talukdar}
\address{Department of Mathematical Sciences, Tezpur University,  Assam 784028, India}
\email{msp21105@tezu.ac.in}

\allowdisplaybreaks
\begin{document}
\begin{abstract}
In the eleventh paper in
the series on MacMahon’s partition analysis, Andrews and Paule \cite{AP07} introduced the $k$-elongated partition diamonds.
Recently, they \cite{AP21} revisited the topic. Let $d_k(n)$ count the partitions obtained by adding the links of the $k$–elongated plane partition diamonds of length $n$. Andrews and Paule \cite{AP21} obtained several generating functions and congruences for $d_1(n)$, $d_2(n)$, and $d_3(n)$. They also posed some conjectures, among which the most difficult one was recently proved by Smoot \cite{Smoot21a}. Da Silva, Hirschhorn, and Sellers \cite{DHS21} further found many congruences modulo certain primes for $d_k(n)$ whereas Li and Yee \cite{LY22} studied the combinatorics of Schmidt type partitions, which can be viewed as partition diamonds. In this article, we give elementary proofs of the remaining conjectures of Andrews and Paule \cite{AP21}, extend some individual congruences found by Andrews and Paule \cite{AP21} and da Silva, Hirschhorn, and Sellers \cite{DHS21} to their respective families as well as find new families of congruences for $d_k(n)$, present a  refinement in an existence result for congruences of $d_k(n)$ found by da Silva, Hirschhorn, and Sellers \cite{DHS21}, and prove some new individual as well as a few families of congruences modulo 5, 7, 8, 11, 13, 16, 17, 19, 23, 25, 32, 49, 64 and  128.

\end{abstract}
\maketitle

\vskip 3mm
\noindent {\footnotesize 2010 Mathematical Reviews Classification Numbers: Primary 11P83; Secondary 05A17.}

\section{Introduction and results}\label{Introduction}
For complex numbers $a$ and $q$ such that $\mid q\mid <1$, we define the infinite $q$-product as
\begin{align*}
\left(a;q\right)_\infty:=\prod_{j=0}^\infty\left(1-aq^j\right).
\end{align*}
For convenience, we set $f_n:=\left(q^n;q^n\right)_\infty$ for integers $n\ge1$.

Andrews and Paule \cite{AP07} in 2007 introduced a combinatorial object called the $k$-elongated partition diamonds. Recently, they \cite{AP21} looked back at  the $k$-elongated partition diamonds. Let $d_k(n)$ count the partitions obtained by adding the links of the $k$–elongated plane partition diamonds of length $n$. Then the generating function for $d_k(n)$ is given by
\begin{align*}
\sum_{n=0}^\infty d_k(n)q^n=\dfrac{f_2^k}{f_1^{3k+1}}.
\end{align*} Andrews and Paule \cite{AP21} found some elegant generating functions for $d_1(n)$, $d_2(n)$, and $d_3(n)$. They also proved many Ramanujan type congruences modulo 2,  3,  4, 5, 8, 9, 27, and  243, mainly by using the Mathematica package \textit{RaduRK} developed by Smoot \cite{Smoot21a}, which uses Radu’s Ramanujan-Kolberg algorithm \cite{Radu15}.

Andrews and Paule \cite{AP21} conjectured some congruences in their paper. For example, they conjectured that for all $n\geq1$ and $k\geq1$ such that $8n\equiv1\pmod{3^k}$, then
\begin{align}
\label{AP conj1}d_{2}\left(n\right)& \equiv0 \pmod{3^k}.
\end{align}
By manipulation of a certain ring of modular functions, Smoot recently \cite{Smoot21b} not only proved but refined \eqref{AP conj1} as
\begin{align*}
d_{2}\left(n\right)& \equiv0 \pmod{3^{\floor*{k/2}+1}}
\end{align*}
for all $n\geq1$ and $k\geq1$ such that $8n\equiv1\pmod{3^k}$. 

Andrews and Paule \cite{AP21} conjectured some more congruences modulo 81, 243, and 729.
\begin{conjecture} \label{AP Conj Theorem 1} 
For all $n\geq0$, 
\begin{align}
\label{AP conj2}d_2\left(81n+44\right)&\equiv0 \pmod{81},\\
\label{AP conj3}d_2\left(81n+j\right)&\equiv0 \pmod{243},\quad \text{where $j\in\{8,35,62,71\}$,}\\
\label{AP conj4}d_2\left(243n+71\right)&\equiv0 \pmod{729}.
\end{align}
\end{conjecture}

\begin{remark}
Andrews and Paule conjectured another congruence \cite[(7.17)]{AP21}, which is in fact contained in \eqref{AP conj3}. 
\end{remark}

Da Silva, Hirschhorn, and Sellers \cite{DHS21} gave elementary proofs for some of the results of Andrews and Paule \cite{AP21} and discovered new individual congruences as well as some infinite families of congruences for $d_{k}(n)$ modulo certain primes. For example, for prime $p\geq5$, let $r$, $1\leq r\leq p-1$, be a quadratic nonresidue modulo $p$. Then for all $n\geq0$ and $N\geq1$,
\begin{align*}
d_{p^N-1}\left(pn+r\right)&\equiv0 \pmod{p^N}.
\end{align*}                                                                                                                                                                                                                                                                                                                                                                                                              Additionally, they \cite{DHS21} proved the following overarching theorem, which generalizes the Ramanujan type congruences modulo prime $p$ with arithmetic progression $p$.
                                                                                                                                                                                                                                                                                                                                                                                                              \begin{theorem}
\label{DHS theorem1}                                                                                                                                                                                                                                                                                                                                                                                                              Let $p$ be a prime, $k\geq1$, $j\geq0$, $N\geq1$ and $r$ be an integer such that \quad \quad \quad \quad \quad $1\leq r\leq p-1$. If, for all $n\geq0$,
\begin{align*}
d_k\left(pn+r\right)&\equiv0 \pmod{p^N},
\end{align*}
then for all $n\geq0$,
\begin{align*}
d_{p^Nj+k}\left(pn+r\right)&\equiv0 \pmod{p^N}.
\end{align*}
\end{theorem}

Andrews and Paule \cite{AP21} considered the partition diamonds as the Schmidt type partitions and in \cite{LY22}, Li and Yee found generating functions for Schmidt $k$-partitions and
unrestricted Schmidt $k$-partitions in unified combinatorial ways.

In this article, in Section \ref{Proof Sect Conj AP}, we give an elementary proof of Conjecture \ref{AP Conj Theorem 1}. In Section \ref{Extnsn of AP and DHS}, we  extend some individual congruences found by Andrews and Paule \cite{AP21} and da Silva, Hirschhorn, and Sellers \cite{DHS21} to their families.  In Section \ref{New Families}, we also find some new families of congruences for $d_k(n)$ modulo 8, 16, 32, 64, and 128. In Section \ref{The refined existence}, we present a  refinement in an existence result for congruences of $d_k(n)$ found by da Silva, Hirschhorn, and Sellers \cite{DHS21}. Finally in Section \ref{The last individuals}, we prove some new individual as well as a few families of congruences modulo 5, 7, 11, 13, 17, 19, 23, 25, and  49.

\section{Proof of Conjecture \ref{AP Conj Theorem 1} of Andrews and Paule}\label{Proof Sect Conj AP} 
First, we recall the following 3-dissections of $f_1^2/f_2$,  $f_2^2/f_1$, $1/f_1^3$, and $f_1f_2$ from \cite[(14.3.2), (14.3.3), (39.2.8), and (14.3.1)]{Hirschhorn17}
\begin{align}
\label{3 dissection f1^2/f2}\frac{f_1^2}{f_2}&=\frac{f_9^2}{f_{18}}-2q \frac{f_3 f_{18}^2}{f_6 f_9},\\
\label{3 dissection f2^2/f1} \frac{f_2^2}{f_1}&=\frac{f_6 f_9^2}{f_3 f_{18}}+q \frac{f_{18}^2}{f_9},\\
\label{3 dissection f1/f2^2}\frac{f_1}{f_2^2}&=\frac{f_3^2 f_9^3}{f_6^6}-q\frac{f_3^3 f_{18}^3}{f_6^7}+q^2\frac{f_3^4 f_{18}^6 }{f_6^8 f_9^3},\\
\label{3 dissection f1^3}f_1^3&=a\left(q^3\right) f_3-3 q f_9^3,\\
\label{3 dissection 1/f1^3}\frac{1}{f_1^3}&= a^2\left(q^3\right)\frac{f_9^3}{f_3^{10}}+3qa^2\left(q^3\right) \frac{f_9^6}{f_3^{11}}+9 q^2\frac{f_9^9}{f_3^{12}},\\
\label{3 dissection f1f2}f_1f_2&=\frac{f_6 f_9^4}{f_3 f_{18}^2}-qf_9 f_{18} -2q^2 \frac{f_3 f_{18}^4}{f_6 f_9^2},
\end{align}
where $a(q)$ is Borweins' cubic theta function defined by $\displaystyle{a(q):=\sum_{j,k=-\infty}^\infty q^{m^2+mn+n^2}}$. 

We have
\begin{align}
\label{Conjecture Step 1} \sum_{n=0}^\infty d_2(n)q^n= \frac{f_2^2}{f_1^7} = \frac{f_2^2}{f_1}\cdot \frac{1}{f_1^6}.
\end{align}
We apply \eqref{3 dissection f2^2/f1} and \eqref{3 dissection 1/f1^3} in \eqref{Conjecture Step 1} to extract the terms involving $q^{3n+2}$ to obtain
\begin{align}
\label{Conjecture Step 2} \sum_{n=0}^\infty d_2(3n+2)q^n= 27a^2(q)\frac{f_2f_3^{14}}{f_1^{23}f_6}+6a^3(q)\frac{f_3^8f_6^{2}}{f_1^{21}}+81q\frac{f_3^{17}f_6^{2}}{f_1^{24}}.
\end{align}

Here, the following identities from \cite[Section 22.10 and (21.3.2)]{Hirschhorn17} again, come into our use
\begin{align}
\label{a(q) to f} a(q)&=\frac{f_1^3}{f_3}+9q \frac{f_9^3}{f_3},\\
\label{a^3(q) to f}a^3(q)&=\frac{f_1^9}{f_3^3}+27 q \frac{f_3^9}{f_1^3},\\
\label{The Main Theta Function Identity}8q \frac{f_1^3 f_6^5}{f_2^3 f_3}&=\dfrac{f_3^8}{f_6^4}-\dfrac{f_1^8}{f_2^4}.
\end{align}
Using \eqref{a(q) to f} and \eqref{a^3(q) to f} in \eqref{Conjecture Step 2}, we have
\begin{align}
\label{Conjecture Step 3}\sum_{n=0}^\infty d_2(3n+2)q^n&= 6\frac{f_3^5f_6^{2}}{f_1^{12}}+27\frac{f_2f_3^{12}}{f_1^{17}f_6}+243q\frac{f_3^{17}f_6^{2}}{f_1^{24}}\notag\\
&\quad +486q\frac{f_2f_3^{12}f_9^3}{f_1^{20}f_6}+2187q^2\frac{f_2f_3^{12}f_9^6}{f_1^{23}f_6}.
\end{align}

 Under modulo $729$, \eqref{Conjecture Step 3} becomes
\begin{align}
\label{Conjecture Step 8}\sum_{n=0}^\infty d_2(3n+2)q^n&\equiv 6\frac{f_3^5f_6^{2}}{f_1^{12}}+27\frac{f_2f_3^{12}}{f_1^{17}f_6}+243q\frac{f_3^{17}f_6^{2}}{f_1^{24}} +486q\frac{f_2f_3^{12}f_9^3}{f_1^{20}f_6}\notag\\
&\equiv 6f_1^{231}\frac{f_6^{2}}{f_3^{76}} +27f_1^{10}f_2\frac{f_3^{3}}{f_6}+243qf_3^9f_6^2\notag\\
&\quad +486qf_1f_2\frac{f_3^5f_9^3}{f_6}\pmod{729}.
\end{align}

 Using \eqref{3 dissection f1^3} and \eqref{3 dissection f1f2}, we find that
 \begin{align}
\label{Conjecture Step 9}  f_1^{231}&\equiv a^{77}\left(q^3\right) f_3^{77}+12  q a^{76}\left(q^3\right) f_3^{76}f_9^3 +90 q^2a^{75}\left(q^3\right) f_3^{75}f_9^6  \notag\\
 &\quad +54q^3 a^{74}\left(q^3\right)f_3^{74} f_9^9  +162 q^4  a^{73}\left(q^3\right)f_3^{73} f_9^{12} \pmod{243},\\
\label{Conjecture Step 10}f_1^{10}f_2&\equiv 18q^3a^2\left(q^3\right)\frac{ f_3^3 f_9f_{18}^4 }{f_6}+9 q^2a^2\left(q^3\right) f_3^2f_9^4 f_{18}  -2q^2a^3\left(q^3\right)\frac{ f_3^4 f_{18}^4 }{f_6 f_9^2}\notag\\
&\quad -9qa^2\left(q^3\right)\frac{ f_3 f_6 f_9^7 }{f_{18}^2}-q a^3\left(q^3\right) f_3^3f_9 f_{18}   +a^3\left(q^3\right)\frac{ f_3^2 f_6 f_9^4}{f_{18}^2}\pmod{27}.
\end{align}
Therefore, invoking \eqref{3 dissection f1f2}, \eqref{Conjecture Step 9}, and \eqref{Conjecture Step 10} in \eqref{Conjecture Step 8}, and then extracting the terms that involve $q^{3n+2}$, we have
\begin{align}
\label{Conjecture Step 11}\sum_{n=0}^\infty d_2(9n+8)q^n&\equiv 540a^{75}(q)\frac{  f_2^2 f_3^6}{f_1}-54a^3(q)\frac{ f_1^7 f_6^4 }{f_2^2 f_3^2}-243a^2(q)\frac{ f_1^5f_3^4 f_6 }{f_2}\notag\\
&\equiv 540a^{75}(q)\frac{  f_2^2 f_3^6}{f_1}-54a^3(q)\frac{ f_1^7 f_6^4 }{f_2^2 f_3^2}-243\frac{ f_1^2f_3^5 f_6 }{f_2}\pmod{729},
\end{align}
where we use the fact that $a(q)\equiv 1 \pmod{3}$, which is clear from \eqref{a(q) to f}.

From \eqref{Conjecture Step 11}, we have
\begin{align*}
\sum_{n=0}^\infty d_2(9n+8)q^n
&\equiv 54a^{75}(q)\frac{  f_2^2 f_3^6}{f_1}-54a^3(q)\frac{ f_1^7 f_6^4 }{f_2^2 f_3^2}\pmod{243},
\end{align*}
which, due to $a^3(q)\equiv 1\pmod{9}$, reduces to
\begin{align*}
\sum_{n=0}^\infty d_2(9n+8)q^n &\equiv 540\frac{f_2^2f_3^{12}}{f_1^{19}}-54\frac{f_1^7f_6^{4}}{f_2^{2}f_3^{2}}\equiv 54\frac{f_2^2f_6^{4}}{f_1f_3^{2}}\left(\frac{f_3^8}{f_6^4}-\frac{f_1^8}{f_2^4}\right)\pmod {243}.
\end{align*}
The above identity, on account of \eqref{The Main Theta Function Identity} and \eqref{3 dissection f1^2/f2}, becomes
\begin{align}
\label{Conjecture Step 5} \sum_{n=0}^\infty d_2(9n+8)q^n &\equiv 432q\frac{f_1^2f_6^{9}}{f_2f_3^3}\equiv 432q\frac{f_6^{9}f_9^{2}}{f_3^3f_{18}}-864q^2\frac{f_6^{8}f_{18}^{2}}{f_3^2f_{9}}\pmod {243}.
\end{align}
Therefore, from the above identity, we evidently have 
\begin{align}
\label{Conjecture Step 6} d_2(27n+8) &\equiv 0 \pmod{243}
\end{align}
for all $n\ge0$.

To prove Theorem \ref{AP Conj Theorem 1}, we require the following 3-dissection of $a(q)$, which was proved by Hirschhorn, Garvan, and Borwein \cite{HGB93}
\begin{align}
\label{3 dissection a(q)} a(q)=a\left(q^3\right)+6q \dfrac{f_9^3}{f_3}.
\end{align}
Using \eqref{3 dissection f2^2/f1}, \eqref{3 dissection f1/f2^2}, \eqref{3 dissection f1^3}, and \eqref{3 dissection a(q)}, we obtain
\begin{align}
\label{Conjecture Step 12} a^{75}(q)\frac{  f_2^2 f_3^6}{f_1}&\equiv 18q^2 a^{74}\left(q^3\right) f_3^5 f_9^2 f_{18}^2 +q a^{75}\left(q^3\right)\frac{ f_3^6 f_{18}^2}{f_9}+18 qa^{74}\left(q^3\right)\frac{ f_3^4 f_6 f_9^5}{f_{18}}\notag\\
&\quad +a^{75}\left(q^3\right)\frac{ f_3^5 f_6 f_9^2}{f_{18}} \pmod{27},\\
\label{Conjecture Step 13}a^3(q)\frac{f_1^7f_6^4}{f_2^2f_3^2}&\equiv 9 q^4a^3\left(q^3\right)\frac{ f_3^2f_9^3 f_{18}^6  }{f_6^4}-9 q^3a^3\left(q^3\right)\frac{ f_3 f_9^6f_{18}^3  }{f_6^3}+12q^3a^4\left(q^3\right)\frac{ f_3^3 f_{18}^6}{f_6^4}\notag\\
&\quad +9q^2a^3\left(q^3\right)\frac{f_9^9 }{f_6^2}-12q^2a^4\left(q^3\right)\frac{ f_3^2f_9^3 f_{18}^3 }{f_6^3}+q^2a^5\left(q^3\right)\frac{f_3^4 f_{18}^6}{f_6^4 f_9^3}\notag\\
&\quad +12qa^4\left(q^3\right)\frac{f_3 f_9^6}{f_6^2}-qa^5\left(q^3\right)\frac{ f_3^3 f_{18}^3}{f_6^3}+a^5\left(q^3\right)\frac{ f_3^2 f_9^3}{f_6^2}\pmod{27}.
\end{align}

Employing \eqref{3 dissection f1^2/f2}, \eqref{Conjecture Step 12}, and  \eqref{Conjecture Step 13} in \eqref{Conjecture Step 11} and then extracting the terms containing $q^{3n+1}$, we find that

\begin{align}
\sum_{n=0}^\infty d_2(27n+17)q^n &\equiv 540 a^{75}(q)\frac{ f_1^6 f_6^2}{f_3}+9720a^{74}(q)\frac{  f_1^4 f_2 f_3^5}{f_6}-648a^4(q)\frac{ f_1 f_3^6 }{f_2^2}\notag\\
&\quad-243a^3(q)\frac{ f_1^6 f_6^2 }{f_3}+54a^5(q)\frac{f_1^3 f_6^3 }{f_2^3}-486qa^3(q)\frac{ f_1^2 f_3^3 f_6^6 }{f_2^4} \notag\\
\label{Conjecture Step 14}  &\equiv 540 a^{75}(q)\frac{ f_1^6 f_6^2}{f_3}+9720\frac{  f_1 f_2 f_3^6}{f_6}-648a(q)\frac{ f_1 f_3^6 }{f_2^2}\notag\\
&\quad-243 f_3 f_6^2 +54a^5(q)\frac{f_1^3 f_6^3 }{f_2^3}-486q\frac{ f_1^2 f_3^3 f_6^5 }{f_2} \pmod{729},
\end{align}
where we use the facts that $a(q)\equiv 1 \pmod{3}$ and $a^3(q)\equiv 1 \pmod{9}$, which is evident from \eqref{a^3(q) to f}.

From \eqref{Conjecture Step 14} and the fact that $a(q)\equiv1\pmod{3}$, we find that
\begin{align}
\label{Conjecture Step New 2} \sum_{n=0}^\infty d_2(27n+17)q^n &\equiv 54\left( a^{75}(q)\frac{ f_1^6 f_6^2}{f_3}+a^5(q)\frac{f_1^3 f_6^3 }{f_2^3}\right)\equiv 27f_3f_6^{2}\pmod {81},
\end{align}
which immediately proves \eqref{AP conj2}.

With the help of \eqref{3 dissection f1/f2^2}, \eqref{3 dissection f1^3},  and \eqref{3 dissection 1/f1^3}, we have
\begin{align*}
a^{75}(q)\frac{ f_1^6 f_6^2}{f_3}&\equiv a^{77}\left(q^3\right) f_3 f_6^2+12q a^{76}\left(q^3\right) f_6^2 f_9^3 + 9 q^2 a^{75}\left(q^3\right)\frac{ f_6^2 f_9^6}{f_3}\pmod{27},\\
a(q)\frac{ f_1 f_3^6 }{f_2^2}&\equiv a\left(q^3\right)\frac{f_3^8 f_9^3 }{f_6^6}+6q\frac{f_3^7 f_9^6  }{f_6^6}-qa\left(q^3\right)\frac{ f_3^9f_{18}^3 }{f_6^7}-6q^2\frac{ f_3^8f_9^3 f_{18}^3  }{f_6^7}\\
&\quad + q^2a\left(q^3\right)\frac{ f_3^{10}f_{18}^6  }{f_6^8 f_9^3}+6 q^3\frac{f_3^9f_{18}^6 }{f_6^8} \pmod{9},\\
a^5(q)\frac{f_1^3 f_6^3 }{f_2^3}&\equiv 
a^7\left(q^3\right)\frac{ f_6^3 f_9^9}{f_3^{10} f_{18}^3}
+6q a^6\left(q^3\right)\frac{ f_6^3 f_9^{12} }{f_3^{11} f_{18}^3}
-6 qa^7\left(q^3\right)\frac{ f_6^2 f_9^6 }{f_3^9}\\
&\quad
+18q^2 a^6\left(q^3\right)\frac{ f_6^2 f_9^9 }{f_3^{10}}
+12q^2 a^7\left(q^3\right)\frac{ f_6f_9^3 f_{18}^3  }{f_3^8}
+18 q^3a^6\left(q^3\right)\frac{ f_6f_9^6 f_{18}^3  }{f_3^9}\\
&\quad
-8q^3 a^7\left(q^3\right)\frac{ f_{18}^6}{f_3^7}
+6q^4 a^6\left(q^3\right)\frac{f_9^3 f_{18}^6  }{f_3^8}
\pmod{27}.
\end{align*}

We now apply the above three identities as well as \eqref{3 dissection f1^2/f2},  \eqref{3 dissection f1/f2^2}, \eqref{3 dissection f1^3},  \eqref{3 dissection 1/f1^3}, and \eqref{3 dissection f1f2} in \eqref{Conjecture Step 14}. Then from the resulting identity, we extract the terms that involve $q^{3n+2}$ to arrive at

\begin{align*}
\sum_{n=0}^\infty d_2(81n+71)q^n &\equiv
-2430\frac{f_1^8 f_3^3 f_6^3 }{f_2^7}
+4860a^{75}(q)\frac{  f_2^2 f_3^6}{f_1}
+648a^7(q)\frac{  f_2 f_3^3 f_6^3}{f_1^8}\notag\\
&\quad
+972 a^6(q)\frac{ f_2^2 f_3^9}{f_1^{10}}
+9720a(q)\frac{ f_1^7 f_6^4}{f_2^2 f_3^2}
-648 a(q)\frac{f_1^{10} f_6^6 }{f_2^8 f_3^3}\notag\\
&\quad
-486 q \frac{ f_1^{11} f_6^{12} }{f_2^{10} f_3^6}
\pmod{729}.
\end{align*}
The above identity, again with the aid of $a(q)\equiv 1 \pmod{3}$, $a^3(q)\equiv 1 \pmod{9}$, and $a(q)\equiv f_1^3/f_3 \pmod{9}$, which follows from \eqref{a(q) to f}, can be rewritten as
\begin{align}
\sum_{n=0}^\infty d_2(81n+71)q^n &\equiv
-2430\frac{f_1^2 f_3^5 f_6 }{f_2}
+5832\frac{  f_2^2 f_3^6}{f_1}
+648a(q)\frac{  f_2 f_3^3 f_6^3}{f_1^8}
+9720\frac{ f_1 f_6^4}{f_2^2}\notag\\
&\quad
-648 a(q)\frac{f_1 f_6^6 }{f_2^8}
-486 q \frac{ f_1^{2} f_6^{9} }{f_2 f_3^3}\notag
\\
&\equiv
-2430\frac{f_1^2 f_3^5 f_6 }{f_2}
+648\frac{  f_2 f_3^2 f_6^3}{f_1^5}
+9720\frac{ f_1 f_6^4}{f_2^2}
-648 \frac{f_2f_3^2 f_6^3 }{f_1^5}
-486 q \frac{ f_1^{2} f_6^{9} }{f_2 f_3^3}
\notag\\
&\equiv
-2430\frac{f_1^2 f_3^5 f_6 }{f_2}
+9720\frac{ f_1 f_6^4}{f_2^2}
-486 q \frac{ f_1^{2} f_6^{9} }{f_2 f_3^3}
\pmod{729}.\notag
\end{align}
From the above identity, we have
\begin{align}
\label{Conjecture Step 7} d_2(81n+71) &\equiv 0 \pmod{243}
\end{align}
for all $n\ge0$. Congruences \eqref{Conjecture Step 6} and \eqref{Conjecture Step 7} together is \eqref{AP conj3}.

Employing \eqref{3 dissection f1^2/f2} and \eqref{3 dissection f1/f2^2} in the identity above, we have
\begin{align*}
\sum_{n=0}^\infty d_2(81n+71)q^n 
&\equiv
-2430\frac{f_3^5 f_6 f_9^2}{f_{18}}
+9720\frac{f_3^2 f_9^3}{f_6^2}
-486 q \frac{f_6^9 f_9^2}{f_3^3 f_{18}}
+4860 q \frac{f_3^6 f_{18}^2}{f_9}\\
&\quad
-9720 q \frac{f_3^3 f_{18}^3}{f_6^3}
+972 q^2 \frac{f_6^8 f_{18}^2}{f_3^2 f_9}
+9720 q^2 \frac{f_3^4 f_{18}^6}{f_6^4 f_9^3}
\pmod{729},
\end{align*}
from which, it follows that
\begin{align*}
\sum_{n=0}^\infty d_2(243n+71)q^n 
&\equiv
-2430\frac{f_1^5 f_2 f_3^2}{f_{6}}
+9720\frac{f_1^2 f_3^3}{f_2^2}\equiv
-2430\frac{f_1^2  f_3^3}{f_{2}^2}
+9720\frac{f_1^2 f_3^3}{f_2^2}\\
&\equiv
7290\frac{f_1^2 f_3^3}{f_2^2}
\pmod{729}.
\end{align*}
The above identity clearly implies \eqref{AP conj4}. \qed

\section{Families for individual congruences of Andrews and Paule \cite{AP21} and da Silva \textit{et. al.} \cite{DHS21}}\label{Extnsn of AP and DHS}
In this section, we extend some of the individual congruences of Andrews and Paule \cite{AP21} and da Silva, Hirschhorn, and Sellers \cite{DHS21} to certain families of congruences in the following theorem.
 
\begin{theorem}\label{Extnsn Thm of AP and DHS} 
For all $n\ge0$ and $j\ge0$,
\begin{align}
\label{AP Indi1} d_{4j+3}\left(4n+2\right)&\equiv 0 \pmod{2},\\
\label{AP Indi2}d_{4j+3}\left(4n+3\right)&\equiv 0 \pmod{4},\\
\label{DHS Indi1}d_{8j+7}\left(4n+2\right)&\equiv 0 \pmod{4},\\
\label{DHS Indi2}d_{8j+7}\left(8n+5\right)&\equiv 0 \pmod{4},\\
\label{DHS Indi3}d_{16j+3}\left(16n+9\right)&\equiv 0 \pmod{4},\\
\label{DHS Indi4}d_{8j+7}\left(4n+3\right)&\equiv 0 \pmod{8},\\
\label{DHS Indi5}d_{32j+7}\left(8n+4\right)&\equiv 0 \pmod{8},\\
\label{DHS Indi6}d_{9j+8}\left(9n+3\right)&\equiv 0 \pmod{9},\\
\label{AP Indi3}d_{27j+2}\left(9n+8\right)&\equiv 0 \pmod{27},\\
\label{AP Indi4}d_{243j+2}\left(27n+8\right)&\equiv 0 \pmod{243}.
\end{align}
\end{theorem}
Note that the individual cases when $j=0$ in \eqref{AP Indi1}, \eqref{AP Indi2}, \eqref{AP Indi3}, \eqref{AP Indi4} and the remainders of the above theorem were proved by Andrews and Paule \cite{AP21} and da Silva, Hirschhorn, and Sellers \cite{DHS21}, respectively.

\begin{proof}
To prove \eqref{AP Indi1}, we first find the following exact generating function
\begin{align}
\label{Exact 1}&\sum_{n=0}^\infty d_{4j+3}(4n+2)q^n \notag \\
&=2 \Bigg\{\dfrac{f_2^{27}f_4^{39}}{f_1^{55}f_8^{18}}\sum_{k=0}^{\floor*{(3j+2)/2}} \sum_{m=0}^{\floor*{(13j+10)/2}}2^{2(2k+m)}\binom{6j+5}{4k}\binom{13j+11}{2m+1}q^{k+m}\dfrac{F_2F_4}{F_1F_8} \notag \\
&\quad + \dfrac{f_2^{13}f_4^{53}}{f_1^{51}f_8^{22}}\sum_{k=0}^{\floor*{(3j+1)/2}} \sum_{m=0}^{\floor*{(13j+11)/2}}2^{2(2k+m)+1}\binom{6j+5}{4k+2}\binom{13j+11}{2m}q^{k+m}\dfrac{F_2F_4}{F_1F_8}\Bigg\},
\end{align}
where $F_1:=f_1^{65j-8k}$, $F_2:=f_2^{30j-24k+4m}$, $F_4:=f_4^{53j+16k-12m}$, and $F_8:=f_8^{26j-8m}$. From \eqref{Exact 1}, \eqref{AP Indi1} is evident.

We need the following 2-dissection from \cite[(1.9.4)]{Hirschhorn17} to establish \eqref{Exact 1}.
\begin{align}
\label{2 dissection of 1/f1^2}\frac{1}{f_1^2}&= \frac{f_8^5}{f_2^5 f_{16}^2}+2q\frac{f_4^2 f_{16}^2}{f_2^5 f_8}.
\end{align}
Using \eqref{2 dissection of 1/f1^2}, we have
\begin{align}
\sum_{n=0}^\infty d_{4j+3}(n)q^n&=\dfrac{f_2^{4j+3}}{f_1^{12j+10}}=f_2^{4j+3}\left(\dfrac{1}{f_1^2}\right)^{6j+5}=f_2^{4j+3}\left(\frac{f_8^5}{f_2^5 f_{16}^2}+2q\frac{f_4^2 f_{16}^2}{f_2^5 f_8}\right)^{6j+5} \notag\\
&=\sum_{k=0}^{6j+5}2^k\binom{6j+5}{k}q^k\dfrac{f_4^{2 k} f_8^{30 j-6 k+25}}{f_2^{26 j+22}f_{16}^{12 j-4 k+10}}\notag\\
&=\dfrac{1}{f_2^{26 j+22}}\Bigg\{\sum_{k=0}^{3j+2}2^{2k}\binom{6j+5}{2k}q^{2k}\dfrac{f_4^{4 k} f_8^{30 j-12 k+25}}{f_{16}^{12 j-8 k+10}}\notag\\
&\quad +\sum_{k=0}^{3j+2}2^{2k+1}\binom{6j+5}{2k+1}q^{2k+1}\dfrac{f_4^{4 k+2} f_8^{30 j-12 k+19}}{f_{16}^{12 j-8 k+6}}\Bigg\}.\notag
\end{align}
Extracting the terms that involve $q^{2n}$ from the above identity, we obtain
\begin{align*}
\sum_{n=0}^\infty d_{4j+3}(2n)q^n&=\left(\dfrac{1}{f_1^{2}}\right)^{13 j+11}\sum_{k=0}^{3j+2}2^{2k}\binom{6j+5}{2k}q^{k}\dfrac{f_2^{4 k} f_4^{30 j-12 k+25}}{f_{8}^{12 j-8 k+10}},
\end{align*}
which again using \eqref{2 dissection of 1/f1^2} can be written as
\begin{align}
\sum_{n=0}^\infty& d_{4j+3}(2n)q^n=\left(\frac{f_8^5}{f_2^5 f_{16}^2}+2q\frac{f_4^2 f_{16}^2}{f_2^5 f_8}\right)^{13j+11}\sum_{k=0}^{3j+2}2^{2k}\binom{6j+5}{2k}q^{k}\dfrac{f_2^{4 k} f_4^{30 j-12 k+25}}{f_{8}^{12 j-8 k+10}}\notag\\
&=\sum_{k=0}^{3j+2}2^{2k}\binom{6j+5}{2k}q^{k}\dfrac{f_2^{4 k} f_4^{30 j-12 k+25}}{f_{8}^{12 j-8 k+10}}\sum_{m=0}^{13j+11}2^m\binom{13j+11}{m}q^m\dfrac{f_4^{2 m} f_8^{65 j-6 m+55}}{f_2^{65 j+55}f_{16}^{26 j-4 m+22}}.\notag
\end{align}
Now, we break the right side of the above identity on the parity of $k$ and $m$ as follows.
\begin{align}
\label{Exact Eve 1} \sum_{n=0}^\infty& d_{4j+3}(2n)q^n\notag\\
&=\Bigg\{\sum_{k=0}^{\floor*{(3j+2)/2}}2^{4k}\binom{6j+5}{4k}q^{2k}\dfrac{f_2^{8 k} f_4^{30 j-24 k+25}}{f_{8}^{12 j-16 k+10}}+\sum_{k=0}^{\floor*{(3j+1)/2}}2^{4k+2}\binom{6j+5}{4k+2}q^{2k+1}\notag\\
&\quad \times\dfrac{f_2^{8 k+4} f_4^{30 j-24 k+13}}{f_{8}^{12 j-16 k+2}}\Bigg\}\Bigg\{\sum_{m=0}^{\floor*{(13j+11)/2}}2^{2m}\binom{13j+11}{2m}q^{2m}\dfrac{f_4^{4 m} f_8^{65 j-12 m+55}}{f_2^{65 j+55}f_{16}^{26 j-8 m+22}}\notag\\
&\quad +\sum_{m=0}^{\floor*{(13j+10)/2}}2^{2m+1}\binom{13j+11}{2m+1}q^{2m+1}\dfrac{f_4^{4 m+2} f_8^{65 j-12 m+49}}{f_2^{65 j+55}f_{16}^{26 j-8 m+18}}\Bigg\}.
\end{align}
From \eqref{Exact Eve 1}, extracting the terms involving $q^{2n+1}$, we deduce \eqref{Exact 1}.

In a similar way, one can find the following generating functions, from which \eqref{AP Indi2}, \eqref{AP Indi3}, and \eqref{DHS Indi4} are evident, respectively.
\begin{align}
\label{Gen 4j+3 4n+3 New} &\sum_{n=0}^\infty d_{4j+3}(4n+3)q^n \notag \\
&=4 \Bigg\{\dfrac{f_2^{7}f_4^{57}}{f_1^{49}f_8^{22}}\sum_{k=0}^{\floor*{(3j+1)/2}} \sum_{m=0}^{\floor*{(13j+11)/2}}2^{2(2k+m)+1}\binom{6j+5}{4k+3}\binom{13j+11}{2m}q^{k+m}\dfrac{F_2F_4}{F_1F_8} \notag \\
&\quad + \dfrac{f_2^{21}f_4^{43}}{f_1^{53}f_8^{18}}\sum_{k=0}^{\floor*{(3j+2)/2}} \sum_{m=0}^{\floor*{(13j+10)/2}}2^{2(2k+m)}\binom{6j+5}{4k+1}\binom{13j+11}{2m+1}q^{k+m}\dfrac{F_2F_4}{F_1F_8}\Bigg\},\\
\label{Gen 8j+7 4n+2 New}&\sum_{n=0}^\infty d_{8j+7}(4n+2)q^n \notag \\
&=4 \dfrac{f_2^{211}}{f_1^{164}f_4^{62}}\Bigg\{\sum_{k=0}^{3j+2} \sum_{m=0}^{\floor*{(13j+12)/2}}2^{4(k+m)}\binom{12j+11}{4k}\binom{13j+12}{2m+1}q^{k+m}\frac{G_2}{G_1G_4} \notag \\
&\quad + \sum_{k=0}^{3j+2} \sum_{m=0}^{\floor*{(13j+11)/2}}2^{4(k+m)}\binom{12j+11}{4k+2}\binom{13j+12}{2m}q^{k+m}\frac{G_2}{G_1 G_4}\Bigg\},\intertext{where $G_1:=f_1^{182 j-8 k-8 m}$, $G_2:=f_2^{242 j-24 k-24 m}$, and $G_4:=f_4^{76 j-16 k-16 m}$,}
\label{Gen 8j+7 4n+3 New}&\sum_{n=0}^\infty d_{8j+7}(4n+3)q^n \notag \\
&=8 \frac{f_2^{205}}{f_1^{162} f_4^{58}}\Bigg\{\sum_{k=0}^{3j+2} \sum_{m=0}^{\floor*{(13j+12)/2}}2^{4(k+m)}\Bigg(\binom{12j+11}{4k+1}\binom{13j+12}{2m+1}+\binom{12j+11}{4k+3}\notag\\
&\quad \times\binom{13j+12}{2m}\Bigg)q^{k+m}\frac{G_2}{G_1 G_4}\Bigg\}.
\end{align}

Note that like the above generating functions,  the exponents of $f_1$ in the generating functions of $d_{8j+7}(4n+1)$,  $d_{16j+3}(4n+1)$, and $d_{32j+7}(4n)$ will also involve $k$. Therefore, the exact generating functions for \eqref{DHS Indi2}, \eqref{DHS Indi3}, and \eqref{DHS Indi5} can not be found as elegantly as the above exact generating functions. So in the following, we give simple proofs for them as well as for the remaining congruences.

The proofs of \eqref{DHS Indi2}, \eqref{DHS Indi3}, and \eqref{DHS Indi5} are similar. So, we prove \eqref{DHS Indi5} only. We have
\begin{align}
\label{32j+7 Eve 1} \sum_{n=0}^\infty d_{32j+7}(n)q^n&=\dfrac{f_2^{32j+7}}{f_1^{96j+22}}\equiv \dfrac{f_1^2}{f_2^{16j+5}}\pmod{8}.
\end{align}
Here, we require the following 2-dissection of $f_1^2$ and $1/f_1^4$ from \cite[(1.9.4) and d (1.10.1)]{Hirschhorn17}.
\begin{align}
\label{2-dissection of f1^2} f_1^2&=\frac{f_2 f_8^5}{f_4^2 f_{16}^2}-2q \frac{f_2 f_{16}^2}{f_8},\\
\label{2-dissection of 1/f1^4}\frac{1}{f_1^4}&=\frac{f_4^{14}}{f_2^{14} f_8^4}+4q \frac{f_4^2 f_8^4}{f_2^{10}}.
\end{align}
Employing \eqref{2-dissection of f1^2} in \eqref{32j+7 Eve 1}, then extracting the terms that involve $q^{2n}$, and using \eqref{2-dissection of 1/f1^4}, we obtain
\begin{align*}
\sum_{n=0}^\infty d_{32j+7}(2n)q^n&\equiv\dfrac{f_4^5}{f_2^{8j+2}f_8^2}\left(\frac{f_4^{14}}{f_2^{14} f_8^4}+4q \frac{f_4^2 f_8^4}{f_2^{10}}\right) \pmod{8},
\end{align*}
which gives
\begin{align*}
\sum_{n=0}^\infty d_{32j+7}(4n)q^n&\equiv\dfrac{f_2^{19}}{f_1^{8j+16}f_4^6}\equiv\dfrac{1}{f_2^{4j-11}f_4^6} \pmod{8}.
\end{align*}
The above identity clearly gives \eqref{DHS Indi5}.

Now, we prove \eqref{DHS Indi6}, \eqref{AP Indi3}, and \eqref{AP Indi4}. Using \eqref{3 dissection f1^2/f2}, we have
\begin{align*}
\sum_{n=0}^\infty d_{9j+8}(n)q^n&=\dfrac{f_2^{9j+8}}{f_1^{27j+25}}\equiv\dfrac{f_6^{3j+3}}{f_3^{9j+9}}\cdot\dfrac{f_1^2}{f_2}\equiv\dfrac{f_6^{3j+3}}{f_3^{9j+9}}\left(\frac{f_9^2}{f_{18}}-2q \frac{f_3 f_{18}^2}{f_6 f_9}\right)\pmod{9},
\end{align*}
which gives
\begin{align*}
\sum_{n=0}^\infty d_{9j+8}(3n)q^n&\equiv\dfrac{f_2^{3j+3}f_3^2}{f_1^{9j+9}f_{6}}\equiv\dfrac{\left(f_2^{3}\right)^{j+1}}{f_3^{3j+1}f_{6}}\pmod{9}.
\end{align*}
Applying \eqref{3 dissection f1^3} in the above identity, then expanding binomially, we find that
\begin{align*}
\sum_{n=0}^\infty d_{9j+8}(3n)q^n&\equiv\dfrac{1}{f_3^{3j+1}f_{6}}\sum_{k=0}^{j+1}(-3)^k\binom{j+1}{k}q^{2k}\left(a\left(q^6\right)f_6\right)^{j-k+1}f_{18}^{3k}\pmod{9}.
\end{align*}
Since in the right side of the above identity, there is no term that involve $q^{3n+1}$, extracting the terms that involve $q^{3n+1}$ from the above identity, we deduce \eqref{DHS Indi6}.

We have
\begin{align*}
\sum_{n=0}^\infty d_{27j+2}(n)q^n&=\dfrac{f_2^{27j+2}}{f_1^{81j+7}}\equiv \dfrac{f_6^{9j}}{f_3^{27j}}\cdot\dfrac{f_2^2}{f_1}\cdot\dfrac{1}{f_1^6}\pmod{27}.
\end{align*}
Using \eqref{3 dissection f2^2/f1} and \eqref{3 dissection 1/f1^3} in the above identity, and then extracting the terms that involve $q^{3n+2}$, we obtain
\begin{align*}
\sum_{n=0}^\infty d_{27j+2}(3n+2)q^n&\equiv\dfrac{f_2^{9j}}{f_1^{27j}}\left(27\frac{f_2 f_3^{12}}{f_1^{17} f_6}+6a^3(q)\frac{f_3^8 f_6^2}{f_1^{21}}+81q \frac{f_3^{17} f_6^2}{f_1^{24}}\right)\\
&\equiv 6\dfrac{f_6^{3j+2}}{f_3^{9j-2}}\cdot\dfrac{1}{f_1^3}\pmod{27}.
\end{align*}
Now, invoking \eqref{3 dissection 1/f1^3} in the above identity, then extracting the terms involving $q^{3n+2}$, we prove \eqref{AP Indi3}.

Finally, we prove \eqref{AP Indi4} using induction on $j$. Andrews and Paule \cite[(7.12)]{AP21} proved \eqref{AP Indi4} for $j=0$. We assume that \eqref{AP Indi4} is true for some integer $j\ge0$. Now,
\begin{align*}
\sum_{n=0}^\infty d_{243(j+1)+2}(n)q^n&=\dfrac{f_2^{243(j+1)+2}}{f_1^{3(243(j+1)+2)+1}}=\dfrac{f_2^{243j+2}}{f_1^{3(243j+2)+1}}\cdot \dfrac{f_2^{243}}{f_1^{729}}=\sum_{n=0}^\infty d_{243j+2}(n)q^n\cdot \dfrac{f_2^{243}}{f_1^{729}}\\
&\equiv \sum_{n=0}^\infty d_{243j+2}(n)q^n\cdot \dfrac{f_6^{81}}{f_3^{243}}\pmod{243}.
\end{align*}
Extracting the terms that involve $q^{3n+2}$ from both sides of the above identity, we have
\begin{align*}
\sum_{n=0}^\infty d_{243(j+1)+2}(3n+2)q^n
&\equiv \sum_{n=0}^\infty d_{243j+2}(3n+2)q^n\cdot \dfrac{f_2^{81}}{f_1^{243}}\pmod{243}.
\end{align*}
Da Silva, Hirschhorn, and Sellers \cite[(21)]{DHS21} showed that $d_{3j+2}(3n+2)\equiv 0\pmod{3}$. Therefore, the above identity can be written as
\begin{align*}
\sum_{n=0}^\infty d_{243(j+1)+2}(3n+2)q^n
&\equiv \sum_{n=0}^\infty d_{243j+2}(3n+2)q^n\cdot \dfrac{f_6^{9}}{f_3^{81}}\pmod{243}.
\end{align*}
Again, extracting the terms that involve $q^{3n+2}$ from both sides of the above identity, we have
\begin{align*}
\sum_{n=0}^\infty d_{243(j+1)+2}(9n+8)q^n
&\equiv \sum_{n=0}^\infty d_{243j+2}(9n+8)q^n\cdot \dfrac{f_2^{9}}{f_1^{81}}\pmod{243}.
\end{align*}
Due to \eqref{AP Indi3}, the above identity is equivalent to
\begin{align*}
\sum_{n=0}^\infty d_{243(j+1)+2}(9n+8)q^n
&\equiv \sum_{n=0}^\infty d_{243j+2}(9n+8)q^n\cdot \dfrac{f_6^{3}}{f_3^{27}}\pmod{243},
\end{align*}
which gives
\begin{align*}
\sum_{n=0}^\infty d_{243(j+1)+2}(27n+8)q^n
&\equiv \sum_{n=0}^\infty d_{243j+2}(27n+8)q^n\cdot \dfrac{f_2^{3}}{f_1^{27}}\pmod{243}.
\end{align*}
Therefore, by the assumption for induction, we see that \eqref{AP Indi4} is true for $j+1$ as well. Thus, \eqref{AP Indi4} is true for all $j\ge0$.
\end{proof}

\section{New families of congruences modulo 8, 16, 32, 64,  and 128}\label{New Families}
In this section, we give new families of congruences taking the advantage of \eqref{Gen 8j+7 4n+2 New} and \eqref{Gen 8j+7 4n+3 New} from Section \ref{Extnsn of AP and DHS}. We also provide an  exact generating function for $d_{16j+15}(4n+3)$ here and use it to prove new congruences.
\begin{theorem}
For all $n\ge0$ and $j\ge0$,
\begin{align}
\label{Cong 8j+7 8n+6} d_{8j+7}\left(8n+6\right)&\equiv 0 \pmod{8},\\
\label{Cong 8j+7 8n+7}d_{8j+7}\left(8n+7\right)&\equiv 0 \pmod{16},\\
\label{Cong 16j+7 8n+6}d_{16j+7}\left(8n+6\right)&\equiv 0 \pmod{16},\\
\label{Cong 16j+7 16n+11}d_{16j+7}\left(16n+11\right)&\equiv 0 \pmod{16},\\
\label{Cong 16j+15 4n+3} d_{16j+15}\left(4n+3\right)&\equiv 0 \pmod{16},\\
\label{Cong 16j+15 8n+6}d_{16j+15}\left(8n+6\right)&\equiv 0 \pmod{16},\\
\label{Cong 16j+15 16n+10}d_{16j+15}\left(16n+10\right)&\equiv 0 \pmod{16},\\
\label{Cong 32j+31 4n+3}d_{32j+31}\left(4n+3\right)&\equiv 0 \pmod{32},\\
\label{Cong 16j+15 8n+7}d_{16j+15}\left(8n+7\right)&\equiv 0 \pmod{64},\\
\label{Cong 32j+31 8n+7}d_{32j+31}\left(8n+7\right)&\equiv 0 \pmod{128}.
\end{align}
\end{theorem}

\begin{proof}
First, we prove \eqref{Cong 16j+15 4n+3}. Similar to \eqref{Gen 4j+3 4n+3 New}--\eqref{Gen 8j+7 4n+3 New}, using \eqref{2 dissection of 1/f1^2} and \eqref{2-dissection of 1/f1^4}, one can find that
\begin{align}
\sum_{n=0}^\infty d_{16j+15}(4n+3)q^n&=8\frac{ f_2^{447}}{f_1^{344} f_4^{134}} \Bigg\{\sum_{k=0}^{6j+5}\sum_{m=0}^{13 j + 12}2^{4 (k+m)}\Bigg(\binom{24 j + 23}{4k+3}\binom{26 j + 25}{2m}\notag\\
&\quad +\binom{24 j + 23}{4k+1}\binom{26 j + 25}{2m+1}\Bigg) q^{k+m}\frac{f_2^{484 j-24 k-24 m}}{f_1^{364 j-8 k-8 m} f_4^{152 j-16 k-16 m}}\Bigg\}.\notag
\end{align}
Now, we separate the right side of the above identity with the cases $(k,m)=(0,0)$ and $(k,m)\neq (0,0)$ as follows.
\begin{align}
\sum_{n=0}^\infty &d_{16j+15}(4n+3)q^n\notag\\
&=8\frac{ f_2^{447}}{f_1^{344} f_4^{134}} \Bigg\{\Bigg(\binom{24 j+23}{3}+\binom{24 j+23}{1} \binom{26 j+25}{1}\Bigg)\frac{ f_2^{484 j}}{f_1^{364 j} f_4^{152 j}}\notag\\
&\quad +\sum_{m=0}^{13 j+11}2^{4 (m+1)}\Bigg(\binom{24 j+23}{3} \binom{26 j+25}{2 m+2}+\binom{24 j+23}{1} \binom{26 j+25}{2 m+3}\Bigg)q^{m+1}\notag\\
&\quad\times \frac{f_2^{484 j-24 m-24}}{f_1^{364 j-8 m-8} f_4^{152 j-16 m-16}}+\sum_{k=0}^{6j+4}2^{4 (k+1)}\Bigg(\binom{24 j+23}{3} \binom{26 j+25}{2 m+2}+\binom{24 j+23}{1}\notag\\
&\quad \times \binom{26 j+25}{2 m+3}\Bigg)q^{k+1}\frac{f_2^{484 j-24 k-24}}{f_1^{364 j-8 k-8} f_4^{152 j-16 k-16}}+\sum_{k=0}^{6j+4}\sum_{m=0}^{13 j + 11}2^{4 (k+m+2)}\Bigg(\binom{24 j+23}{4 k+7}\notag\\
&\quad \times \binom{26 j+25}{2 m+2} +\binom{24 j+23}{4 k+5} \binom{26 j+25}{2 m+3}\Bigg) q^{k+m+2}\frac{f_2^{484 j-24 k-24 m-48}}{f_1^{364 j-8 k-8 m-16} f_4^{152 j-16 k-16 m-32}}\Bigg\}.\notag
\end{align}
On simplifying the above identity, we find that
\begin{align}
\label{Gen 16j+15 4n+3} \sum_{n=0}^\infty &d_{16j+15}(4n+3)q^n\notag\\
&=16\Bigg\{3  (24 j+23) (16 j+17)(j+1)\frac{ f_2^{484 j+447}}{f_1^{364 j+344} f_4^{152 j+134}}+8\sum_{m=0}^{13 j+11}2^{4m}\Bigg(\binom{24 j+23}{3}\notag\\
&\quad \times \binom{26 j+25}{2 m+2}+(24 j+23) \binom{26 j+25}{2 m+3}\Bigg)q^{m+1} \frac{f_2^{484 j-24 m+423}}{f_1^{364 j-8 m+336} f_4^{152 j-16 m+118}}\notag\\
&\quad+8\sum_{k=0}^{6j+4}2^{4k}\Bigg(\binom{24 j+23}{3} \binom{26 j+25}{2 m+2}+(24 j+23)\binom{26 j+25}{2 m+3}\Bigg)q^{k+1}\notag\\
&\quad \times \frac{f_2^{484 j-24 k+423}}{f_1^{364 j-8 k+336} f_4^{152 j-16 k+118}}+128\sum_{k=0}^{6j+4}\sum_{m=0}^{13 j + 11}2^{4 (k+m)}\Bigg(\binom{24 j+23}{4 k+7}\binom{26 j+25}{2 m+2}\notag\\
&\quad  +\binom{24 j+23}{4 k+5} \binom{26 j+25}{2 m+3}\Bigg) q^{k+m+2}\frac{f_2^{484 j-24 k-24 m+399}}{f_1^{364 j-8 k-8 m+328} f_4^{152 j-16 k-16 m+102}}\Bigg\}.
\end{align}
Note that \eqref{Cong 16j+15 4n+3} is evident from \eqref{Gen 16j+15 4n+3}.

Now, we prove \eqref{Cong 32j+31 8n+7}. Replacing $j$ by $2j+1$ and taking modulo $128$ in \eqref{Gen 16j+15 4n+3}, we obtain
\begin{align*}
\sum_{n=0}^\infty d_{32j+31}(4n+3)q^n
&\equiv 96  (48 j+47) (32 j+33)(j+1)\frac{ f_2^{968 j+931}}{f_1^{728 j+708} f_4^{304 j+286}}\\
&\equiv 96  (48 j+47) (32 j+33)(j+1)\frac{ f_2^{604 j+577}}{ f_4^{304 j+286}} \pmod{128}.
\end{align*}
From the above identity, we clearly have \eqref{Cong 32j+31 8n+7}.

Similar to the proof of \eqref{Cong 32j+31 8n+7}, we can obtan \eqref{Cong 8j+7 8n+6}, \eqref{Cong 16j+7 8n+6}, \eqref{Cong 16j+15 8n+6}, and \eqref{Cong 16j+15 16n+10} from \eqref{Gen 8j+7 4n+2 New}, \eqref{Cong 8j+7 8n+7} and \eqref{Cong 16j+7 16n+11} from \eqref{Gen 8j+7 4n+3 New}, and \eqref{Cong 32j+31 4n+3} and \eqref{Cong 16j+15 8n+7} from \eqref{Gen 16j+15 4n+3}.
\end{proof}

\section{An existence result for infinite families of congruences}\label{The refined existence}
In this section, we provide the following theorem that refines Theorem \ref{DHS theorem1}, which was found by da Silva, Hirschhorn, and Sellers \cite{DHS21}. 
\begin{theorem}\label{The refined Theorem} 
Let $p$ be a prime, $k\ge1$, $j\ge0$, $N\ge1$, $M\ge1$, and $r$ be integers such that $1\le r \le p^M-1$. If for all $n\ge0$,
\begin{align*}
d_k\left(p^Mn+r\right)\equiv 0 \pmod{p^N},
\end{align*}
then for all $n\ge0$,
\begin{align*}
d_{p^{M+N-1}j+k}\left(p^Mn+r\right)\equiv 0 \pmod{p^N}.
\end{align*}
\end{theorem}
\begin{proof}
Without loss of generality, we may assume that $\displaystyle{r=\sum_{j=0}^{M-1}p^jr_j}$ for $0\le r_j\le p-1$, because  $\displaystyle{\sum_{j=0}^{M-1}p^jr_j}$ can take any value between 1 and $p^M-1$. For integers  $M\ge1$ (sufficiently large) and $N\ge1$, we have
\begin{align*}
\sum_{n=0}^\infty d_{p^{M+N-1}j+k}(n)q^n &= \dfrac{f_2^{p^{M+N-1}j+k}}{f_1^{3p^{M+N-1}j+3k+1}}\equiv \dfrac{f_{2p}^{p^{M+N-2}j}}{f_{p}^{3p^{M+N-2}j}}\sum_{n=0}^\infty d_k(n)q^n \pmod{p^N}.
\end{align*}
Extracting the terms that involve $q^{pn+r_0}$ from the above identity, we obtain 
\begin{align*}
\sum_{n=0}^\infty d_{p^{M+N-1}j+k}(pn+r_0)q^n 
&\equiv \dfrac{f_{2}^{p^{M+N-2}j}}{f_{1}^{3p^{M+N-2}j}}\sum_{n=0}^\infty d_k(pn+r_0)q^n \\
&\equiv \dfrac{f_{2p}^{p^{M+N-3}j}}{f_{p}^{3p^{M+N-3}j}}\sum_{n=0}^\infty d_k(pn+r_0)q^n \pmod{p^N}.
\end{align*}
Now, extracting the terms that involve $q^{pn+r_1}$ from the above identity, we find that
 \begin{align*}
\sum_{n=0}^\infty d_{p^{M+N-1}j+k}(p^{2}n+r_0+ pr_1)q^n
&\equiv \dfrac{f_{2}^{p^{M+N-3}j}}{f_{1}^{3p^{M+N-3}j}}\sum_{n=0}^\infty d_k(p^{2}n+r_0+ pr_1)q^n \\
&\equiv \dfrac{f_{2p}^{p^{M+N-4}j}}{f_{p}^{3p^{M+N-4}j}}\sum_{n=0}^\infty d_k(p^{2}n+r_0+ pr_1)q^n \pmod{p^N}.
\end{align*}
From the above identity, we extract the terms that contain $q^{pn+r_2}$, and from the resulting identity, we again  extract the terms that contain $q^{pn+r_3}$. It can be seen that after the $M$-th extraction using this iterative scheme, we arrive at
\begin{align*}
\sum_{n=0}^\infty &d_{p^{M+N-1}j+k}(p^{M}n+r_0+ pr_1+\cdots+p^{M-1}r_{M-1})q^n\\
&\equiv \dfrac{f_{2}^{p^{N-1}j}}{f_{1}^{3p^{N-1}j}}\sum_{n=0}^\infty d_k(p^{2}n+r_0+ pr_1+\cdots+p^{M-1}r_{M-1})q^n \pmod{p^N}.
\end{align*}
Therefore, if we assume that $d_k(p^{M}n+r_0+ pr_1+\cdots+p^{M-1}r_{M-1})=d_k(p^{M}n+r)\equiv 0 \pmod{p^N}$, from the above identity, we evidently have

$$d_{p^{M+N-1}j+k}(p^{M}n+r)\equiv 0 \pmod{p^N}.$$
Thus, we complete the proof Theorem \ref{The refined Theorem}.
\end{proof}

\begin{remark}
Theorem \ref{The refined Theorem} is a refinement of Theorem \ref{DHS theorem1} in the sense that it extends individual congruences with arithmetic progressions $p^Mn+r$, $M\ge1$ to their respective families, whereas Theorem \ref{DHS theorem1} extends individual congruences with arithmetic progressions $pn+r$ only to their respective families. For example, for all $n$, we have
\begin{align}
\label{Example help 1} d_7(2n+1)\not\equiv d_7(8n+7)\equiv 0 \pmod{16}.
\end{align}
So, Theorem \ref{DHS theorem1} does not provide any information regarding its extension to an infinite family, whereas Theorem \ref{The refined Theorem} and \eqref{Example help 1} imply that
\begin{align*}
d_{64j+7}(8n+7)\equiv 0 \pmod{16}.
\end{align*}
\end{remark}

\section{New individual and families of congruences modulo 5, 7, 11, 13, 17, 19, 23, 25, and 49}\label{The last individuals} 

In this section, we present some new individual as well as a few families of congruences. Here, we use modular identities of the Rogers-Ramanujan continued fraction, which is defined as
\begin{align*}
R(q):=\dfrac{1}{1} \cplus  \frac{q}{1}\cplus\frac{q^2}{1}\cplus  \frac{q^3}{1}\cplus\contdots
\end{align*}
 a 7-dissection of $f_1$, series representations of certain $q$-products, and an algorithm developed by Radu \cite{Radu09} to prove the following congruences.
\begin{theorem}\label{The Main Theorem} 
For all $n\ge0$ and $j\ge0$, we have
\begin{align}
\label{d1 mod5 1} d_1\left(25n+23\right)&\equiv 0 \pmod{5},\\
\label{d1 mod25 1}d_1\left(125n+j\right)&\equiv 0 \pmod{25},\quad \text{where $j \in \{23,123\}$,}\\
\label{d2 mod5 1}d_2\left(125n+j\right)&\equiv 0 \pmod{5},\quad \text{where $j \in \{97,122\}$,}\\
\label{d1 mod7 1}d_1\left(49n+j\right)&\equiv 0 \pmod{7},\quad \text{where $j \in \{17,31,38,45\}$,}\\
\label{d2 mod7 1}d_2\left(49n+43\right)&\equiv 0 \pmod{7},\\
\label{d3 mod49 1}d_3\left(49n+41\right)&\equiv 0 \pmod{7},\\
\label{d3 mod7 1}d_3\left(343n+j\right)&\equiv 0 \pmod{49},\quad \text{where $j \in \{90,188,237\}$,}\\
\label{d4 mod7 1}d_4\left(343n+j\right)&\equiv 0 \pmod{7},\quad \text{where $j \in \{39,235,284\}$,}\\
\label{d4 mod11 1}d_4\left(121n+96\right)&\equiv 0 \pmod{11},\\
\label{d5 mod11 1}d_5\left(121n+91\right)&\equiv 0 \pmod{11},\\
\label{d7 mod11 1}d_7\left(121n+81\right)&\equiv 0 \pmod{11},\\
\label{d13j+3 mod13 1}d_{13j+3}\left(13n+11\right)&\equiv 0 \pmod{13},\\
\label{d17j+5 mod17 1}d_{17j+5}\left(17n+13\right)&\equiv 0 \pmod{17},\\
\label{d6 mod17 1}d_{6}\left(289n+j\right)&\equiv 0 \pmod{17},\quad \text{where $j \in \{52,69,137,171,$}\notag\\
&\hspace*{4.2cm} \text{$188,205,222,239,273\}$,}\\
\label{d19j+3 mod19 1}d_{19j+3}\left(19n+16\right)&\equiv 0 \pmod{19},\\
\label{d19j+6 mod19 1}d_{19j+6}\left(19n+9\right)&\equiv 0 \pmod{19},\\
\label{d19j+7 mod19 1}d_{19j+7}\left(19n+13\right)&\equiv 0 \pmod{19},\\
\label{d23j+8 mod23 1}d_{23j+8}\left(23n+9\right)&\equiv 0 \pmod{23}.
\end{align}
\end{theorem}

\subsection{Required lemmas}\label{Section Lemma}
Here, we present some background material on the method of Radu \cite{Radu09}. For integers $x$, let $[x]_m$ denote the residue class of $x$ in $\mathbb{Z}/m\mathbb{Z}$, $\mathbb{Z}_m^{*}$ be the set of all invertible elements in $\mathbb{Z}_m$, $\mathbb{S}_m$ denote the set of all squares in $\mathbb{Z}_m^{*}$, and for integers $N\ge1$, we assume that
\begin{align*}
\Gamma&:=\left\{
\begin{pmatrix}
a & b\\
c & d
\end{pmatrix} : a, b, c, d \in \mathbb{Z}, \text{~and~}ad-bc=1
\right\},\\
\Gamma_{\infty}&:=\left\{
\begin{pmatrix}
1 & n\\
0 & 1
\end{pmatrix}  : n\in  \mathbb{Z}
\right\},\\
\Gamma_0(N)&:=\left\{
\begin{pmatrix}
a & b\\
c & d
\end{pmatrix}  \in \Gamma : c\equiv0\pmod{N}
\right\},\\
[\Gamma : \Gamma_0(N)]&:=N\prod_{\ell\mid N} \left( 1+\dfrac{1}{\ell}\right),
\end{align*}
where $\ell$ is a prime.

For integers $M\ge1$, suppose that $R(M)$ is the set of all the integer sequences $(r_\delta):=\left(r_{\delta_1},r_{\delta_2},r_{\delta_3},\ldots,r_{\delta_k}\right)$
indexed by all the positive divisors $\delta$ of $M$, where $1=\delta_1<\delta_2<\cdots<\delta_k=M$. For integers $m\ge1$, $(r_\delta)\in R(M)$, and $t\in\{0,1,2,\ldots,m-1\}$, we define the set $P(t)$ as
\begin{align}
\label{P(t)} P(t):=&\bigg\{t^\prime \in \{0,1,2,\ldots,m-1\} : t^{\prime}\equiv ts+\dfrac{s-1}{24}\sum_{\delta\mid M}\delta r_\delta \pmod{m} \notag\\
&\quad \text{~for some~} [s]_{24m}\in \mathbb{S}_{24m}\bigg\}.
\end{align}
For integers $N\ge1$, $\gamma:=
\begin{pmatrix}
	a  &  b \\
	c  &  d
\end{pmatrix} \in \Gamma$, $(r_\delta)\in R(M)$, and $(r_\delta^\prime)\in R(N)$, we also define
	\begin{align*}
	p(\gamma)&:=\min_{\lambda\in\{0, 1, \ldots, m-1\}}\dfrac{1}{24}\sum_{\delta\mid M}r_{\delta}\dfrac{(\delta (a+ k\lambda c), mc)^2}{\delta m},\\
	p^\prime(\gamma)&:=\dfrac{1}{24}\sum_{\delta\mid N}r_{\delta}^\prime\dfrac{(\delta, c)^2}{\delta}.
	\end{align*}

For integers $m\ge1$; $2\nmid m$, $M\ge1$, $N\ge1$,  $t\in \{0,1,2,\ldots,m-1\}$, $k:=\left(m^2-1,24\right)$, and $(r_{\delta})\in R(M)$, define $\Delta^{*}$ to be the set of all tuples $(m, M, N, t, (r_{\delta}))$ such that all of the following conditions are satisfied
	
	\begin{enumerate}
		\item[1.] Prime divisors of $m$ are also prime divisors of $N$;
		\item[2.] If $\delta\mid M$, then $\delta\mid mN$ for all $\delta\geq1$ with $r_{\delta} \neq 0$;
		\item[3.] $\displaystyle{24\mid kN\sum_{\delta\mid M}\dfrac{r_{\delta} mN}{\delta}}$;
		\item[4.] $\displaystyle{8\mid kN\sum_{\delta\mid M}r_{\delta}}$;
		\item[5.]  $\dfrac{24m}{\left(-24kt-k{\displaystyle{\sum_{\delta\mid M}}{\delta r_{\delta}}},24m\right)} \mid N$.
	\end{enumerate}

The following lemma supports Lemma \ref{Lemma Radu 1} in the proof of Theorem \ref{The Main Theorem}.
\begin{lemma}\label{Lemma Wang 1}\textup{\cite[Lemma 4.3]{Wang17}} Let $N$ or $\frac{1}{2}N$ be a square-free integer, then we have
		\begin{align*}
		\bigcup_{\delta\mid N}\Gamma_0(N)\begin{pmatrix}
		1  &  0 \\
		\delta  &  1
		\end{pmatrix}\Gamma_ {\infty}=\Gamma.
		\end{align*}
	\end{lemma}

We end this section by stating a result of Radu \cite{Radu09}, which is especially useful in completing the proof of Theorem \ref{The Main Theorem} in the final section.
\begin{lemma}\label{Lemma Radu 1}
\textup{\cite[Lemma 4.5]{Radu09}} Suppose that $(m, M, N, t, (r_{\delta}))\in\Delta^{*}$, $(r'_{\delta}):=(r'_{\delta})_{\delta \mid N}\in R(N)$, $\{\gamma_1,\gamma_2, \ldots, \gamma_n\}\subseteq \Gamma$ is a complete set of representatives of the double cosets of $\Gamma_{0}(N) \backslash \Gamma/ \Gamma_\infty$, $\displaystyle{t_{\min}:=\min_{t^\prime \in P(t)} t^\prime}$,
\begin{align}
	\label{Nu} \nu:= \dfrac{1}{24}\left( \left( \sum_{\delta\mid M}r_{\delta}+\sum_{\delta\mid N}r_{\delta}^\prime\right)[\Gamma:\Gamma_{0}(N)] -\sum_{\delta\mid N} \delta r_{\delta}^\prime-\frac{1}{m}\sum_{\delta|M}\delta r_{\delta}\right)
	- \frac{ t_{min}}{m},
 	\end{align}
$p(\gamma_j)+p^\prime(\gamma_j) \geq 0$ for all $1 \leq j \leq n$, and $\displaystyle{\sum_{n=0}^{\infty}A(n)q^n:=\prod_{\delta\mid M}f_\delta^{r_\delta}.}$ If for some integers $u\ge1$, all $t^\prime \in P(t)$, and $0\leq n\leq \lfloor\nu\rfloor$, $A(mn+t^\prime)\equiv0\pmod u$  is true,  then for integers $n\geq0$ and all $t^\prime\in P(t)$, we have $A(mn+t^\prime)\equiv0\pmod u$.
\end{lemma}

\subsection{Proof of Theorem \ref{The Main Theorem}}\label{Proofs Theorem} 

\begin{proof}[Proof of \eqref{d1 mod5 1}]
First, for integers $\alpha\ge 0$ and $\beta$, we let
\begin{align}
\label{P alpha beta} P_{\alpha,\beta}&:=\dfrac{1}{ R(q)^{\alpha+2\beta}R\left(q^2\right)^{2\alpha-\beta}}+(-1)^{\alpha+\beta}q^{2\alpha} R(q)^{\alpha+2\beta}R\left(q^2\right)^{2\alpha-\beta}
\end{align}
and recall two 5-dissections from \cite[pp. 161--165]{Berndt06} as follows
\begin{align}
\label{5-dissection of f1}f_1=&f_{25}\left(\dfrac{1}{R\left(q^5\right)}-q-q^2 R\left(q^5\right)\right),\\
\label{5-dissection of 1/f1}\dfrac{1}{f_1}=&\dfrac{f_{25}^5}{f_{5}^6}\Bigg(\dfrac{1}{R^4\left(q^5\right)}+ \dfrac{q}{R^3\left(q^5\right)}+\dfrac{2q^2}{R^2\left(q^5\right)}+\dfrac{3q^3}{R\left(q^5\right)}+5q^4\notag\\
&-3q^5 R\left(q^5\right)+2q^6 R^2\left(q^5\right)-q^7 R^3\left(q^5\right)+q^8 R^4\left(q^5\right)\Bigg). 
\end{align}

Now,
\begin{align*}
\sum_{n=0}^\infty d_1(n)q^n=\dfrac{f_2}{f_1^4}.
\end{align*}
Employing the 5-dissections of $f_2$ and $1/f_1$ from \eqref{5-dissection of f1} and \eqref{5-dissection of 1/f1} in the above identity, then extracting the terms that involve $q^{5n+3}$, and finally with the help of \eqref{P alpha beta}, we obtain
\begin{align}
\label{RR Step 1} \sum_{n=0}^\infty d_1(5n+3)q^n=\dfrac{f_5^{20}f_{10}}{f_1^{24}}&\Big(-4 P_{3,6}+40  P_{3,5}-105q P_{2,5}-418qP_{2,4}\notag\\
&\quad+1100q  P_{2,3}-1400q^2 P_{1,3}-1840q^2 P_{1,2}\notag\\
&\quad+1200q^2 P_{1,1}-1500q^3P_{0,1}-1015q^3\Big).
\end{align}

From \cite[Lemma 1.3]{BB18} and \cite[(7.4.9)]{Berndt06}, we have
\begin{align}
\label{Modular Relation 1} P_{0,1}&=4q\frac{f_1 f_{10}^5}{f_2 f_5^5},\\
\label{Modular Relation 2}P_{1,1}&=\frac{f_2 f_5^5}{f_1 f_{10}^5}+ 2 q+ 4q^2\frac{f_1 f_{10}^5}{f_2 f_5^5}\\
\label{Modular Relation 3}P_{1,2}&=\frac{f_1^6}{f_5^6}+11 q.
\end{align}
and the following relations hold
\begin{align*}
P_{1,3}&= P_{0,1}P_{1,2}+P_{1,1}, & 
P_{2,3}&= P_{1,1}P_{1,2}-q^2P_{0,1},\\
P_{2,4}&= P_{1,2}^2+2q^2, & 
P_{2,5}&= P_{0,1}P_{2,4}-P_{2,3},\\
P_{3,5}&= P_{1,1}P_{2,4}-q^2P_{1,3}, & 
P_{3,6}&= P_{1,2}P_{2,4}+q^2P_{1,2}.
\end{align*}

Employing \eqref{Modular Relation 1}--\eqref{Modular Relation 3} and the above relations in \eqref{RR Step 1}, we find that
\begin{align*}
\sum_{n=0}^\infty d_1(5n+3)q^n=&40\frac{ f_2 f_5^{13}}{f_1^{13} f_{10}^4}-4\frac{ f_{10} f_5^2}{f_1^6}-470q\frac{ f_{10} f_5^8 }{f_1^{12}}+1875q\frac{ f_2 f_5^{19} }{f_1^{19} f_{10}^4}\notag\\
&+15625q^2\frac{ f_2 f_5^{25} }{f_1^{25} f_{10}^4}-8750q^2\frac{ f_{10} f_5^{14} }{f_1^{18}}-260q^2\frac{ f_{10}^6 f_5^3 }{f_1^{11} f_2}\notag\\
&-7500q^3\frac{ f_{10}^6 f_5^9 }{f_1^{17} f_2}-46875q^3\frac{ f_{10} f_5^{20} }{f_1^{24}}-62500q^4\frac{ f_{10}^6 f_5^{15} }{f_1^{23} f_2},
\end{align*}
which under modulo $5$ gives
\begin{align*}
\sum_{n=0}^\infty d_1(5n+3)q^n\equiv&-4\frac{ f_{10} f_5^2}{f_1^6}\pmod5\\
\equiv&\frac{ f_{10} f_5}{f_1}\pmod5.
\end{align*}
Invoking the $5$-dissection of $1/f_1$ given by \eqref{5-dissection of 1/f1} in the above identity and then extracting the terms involving $q^{5n+4}$, we obtain \eqref{d1 mod5 1}.
\end{proof}

\begin{proof}[Proof of \eqref{d2 mod7 1}]
We have
\begin{align}
\label{Mod 7 Elementary Step 1} \sum_{n=0}^\infty d_2(n)q^n=\dfrac{f_2^2}{f_1^7}\equiv \dfrac{f_2^2}{f_7}\pmod 7.
\end{align}
From \cite[(10.5.1)]{Hirschhorn17}, we recall the following 7-dissection of $f_1$.
\begin{align*}
f_1=&f_{49}\Bigg(\dfrac{\left(q^{14};q^{49}\right)_\infty\left(q^{35};q^{49}\right)_\infty}{\left(q^{7};q^{49}\right)_\infty\left(q^{42};q^{49}\right)_\infty}-q \dfrac{\left(q^{21};q^{49}\right)_\infty\left(q^{28};q^{49}\right)_\infty}{\left(q^{14};q^{49}\right)_\infty\left(q^{35};q^{49}\right)_\infty}\\
&-q^2+q^5\dfrac{\left(q^{7};q^{49}\right)_\infty\left(q^{42};q^{49}\right)_\infty}{\left(q^{21};q^{49}\right)_\infty\left(q^{28};q^{49}\right)_\infty}\Bigg).
\end{align*}
With the help of the above identity, we use the 7-dissection of $f_2^2$ in \eqref{Mod 7 Elementary Step 1} and then extract the terms involving $q^{7n+1}$. This gives
\begin{align}
\label{Mod 7 Elementary Step 2} \sum_{n=0}^\infty d_2(7n+1)q^n\equiv \dfrac{f_{14}^2}{f_1}\pmod 7.
\end{align}

Now, if $p(n)$ counts the unrestricted partitions of an integer $n\ge0$, we have 
\begin{align*}
\sum_{n=0}^\infty p(n)q^n=\dfrac{1}{f_1}
\end{align*}
and one of Ramanujan's famous three partition congruences
\begin{align*}
p(7n+6)\equiv 0 \pmod 7
\end{align*}
for all $n\ge0$.

Therefore, it becomes evident from \eqref{Mod 7 Elementary Step 2} that \eqref{d2 mod7 1} is true.
\end{proof}

\begin{proof}[Proofs of  \eqref{d13j+3 mod13 1}, \eqref{d17j+5 mod17 1}, and \eqref{d19j+6 mod19 1}--\eqref{d23j+8 mod23 1}]
We have
\begin{align*}
\sum_{n=0}^\infty d_3(n)q^n &= \frac{f_2^3}{f_1^{10}}\equiv  \frac{f_1^3f_2^3}{f_{13}}\pmod{13}.
\end{align*}
Using Jacobi's famous identity \cite[(1.3.24)]{Berndt06} in the above identity,
\begin{align}
\label{d3 ele step1}\sum_{n=0}^\infty d_3(n)q^n &\equiv \frac{1}{f_{13}}\sum_{j=0}^\infty (-1)^j (2j+1)q^{j(j+1)/2}\sum_{k=0}^\infty (-1)^k (2k+1)q^{k(k+1)}\notag\\
&\equiv \frac{1}{f_{13}}\sum_{j,k=0}^\infty (-1)^{j+k} (2j+1)(2k+1)q^{j(j+1)/2 + k(k+1)} \pmod{13}.
\end{align}

Now, 
\begin{align*}
8\left( \frac{j(j+1)}{2} + k(k+1)\right)+3
&= (2j+1)^2 +2(2k+1)^2.
\end{align*}
If $j(j+1)/2+ k(k+1) = 13n+11$ for some integer $n\geq0$, the above equality gives 
\begin{align*}
(2j+1)^2 +2(2k+1)^2 \equiv0 \pmod{13}.
\end{align*}
 Therefore, $2j+1 \equiv0 \pmod{13}$ and $2k+1 \equiv0 \pmod{13}$. Otherwise, we have \linebreak $(2j+1)^2 \equiv1,3,4,9,10,12 \pmod{13}$, which gives $(2j+1)^2 +2(2k+1)^2 \not\equiv 0 \pmod{13}.$ This is a contradiction.

Finally, extracting the terms that involve $q^{13n+11}$ from \eqref{d3 ele step1}, we find that for all $n\geq0$, 
\begin{align}
\label{d3 ele step final}d_3(13n+11) \equiv 0 \pmod {13}.
\end{align}

Theorem \ref{DHS theorem1} and \eqref{d3 ele step final} ensure \eqref{d13j+3 mod13 1}.

Congruences \eqref{d17j+5 mod17 1} and \eqref{d19j+6 mod19 1}--\eqref{d23j+8 mod23 1} can be proved similarly as above. So, we do not go in detail but provide  the following product-to-sum identities and a chart useful for their proofs.
\end{proof}
\begin{align}
\label{jacobi1}f_1^{3}&=\sum_{j=0}^\infty (-1)^j (2j+1)q^{j(j+1)/2}, \\
\label{f2^2 by f1}\dfrac{f_2^2}{f_1}&=\sum_{j=0}^\infty q^{j(j+1)/2},\\
\label{f2^5 by f1^2}\dfrac{f_2^5}{f_1^2}&= \sum_{j=-\infty}^\infty (-1)^j(3j+1)q^{3j^2+2j}.
\end{align}
\begin{longtable}{cc}
\hline
 Congruence  & Used product-to-sum identities
 \T\B\\  \hline
\eqref{d17j+5 mod17 1}  &\eqref{jacobi1}, \eqref{f2^5 by f1^2} \T\B\ \\
\eqref{d19j+6 mod19 1} &\eqref{jacobi1} \T\B\ \\
\eqref{d19j+7 mod19 1} &\eqref{f2^2 by f1}, \eqref{f2^5 by f1^2} \T\B\ \\ 
\eqref{d23j+8 mod23 1} &\eqref{jacobi1}, \eqref{f2^5 by f1^2} \T\B\ \\
\hline
\end{longtable}

\begin{proof}[Proofs of the remaining congruences of Theorem \ref{The Main Theorem}]
Proofs of \eqref{d1 mod25 1}--\eqref{d1 mod7 1}, \eqref{d3 mod49 1}--\eqref{d7 mod11 1}, \eqref{d6 mod17 1}, and \eqref{d19j+3 mod19 1} are similar. We  elaborate the proof of \eqref{d1 mod25 1}  only. We have
\begin{align}
\label{Modular Proof Step 1} \sum_{n=0}^\infty d_1(n)q^n=\dfrac{f_2}{f_1^4}\equiv \dfrac{f_1^{21}f_2}{f_1^{25}}\equiv \dfrac{f_1^{21}f_2}{f_5^{5}}\pmod{25}.
\end{align}

Using Conditions 1--5, it is clear that $(m,M,N,t,(r_\delta))=(125,10,10,23,(21,1,-5,0))$ $\in \Delta^{*}$. So, by \eqref{P(t)}, we have $P(t)=\{23,123\} $. Lemma \ref{Lemma Wang 1} gives that
$\left\{
\begin{pmatrix}
1 & 0\\
\delta & 1
\end{pmatrix} : \delta\mid N
\right\}$ is a complete set of representatives of the double cosets in $\Gamma_{0}(N) \backslash \Gamma/ \Gamma_\infty$. Using $(r_\delta^\prime)=(18,0,0,0)$, \eqref{Nu}, and \textit{Mathematica}, we find that
\begin{align*}
p\left(\begin{pmatrix}
1 & 0\\
\delta & 1
\end{pmatrix}\right)+p^{\prime}\left(\begin{pmatrix}
1 & 0\\
\delta & 1
\end{pmatrix}\right)&\ge 0\quad \text{for all $\delta\mid N$,}\\
\floor*{\nu}&=25,\\
d_1(125n+j)&\equiv0 \pmod{25}
\end{align*}
for $j\in\{23,123\}$ are true for all $0\le n \le \floor*{\nu}$. Therefore, by Lemma \ref{Lemma Radu 1} and \eqref{Modular Proof Step 1}, \eqref{d1 mod25 1} is true. The proofs of \eqref{d2 mod5 1},  \eqref{d1 mod7 1} \eqref{d3 mod49 1}--\eqref{d7 mod11 1}, \eqref{d6 mod17 1}, and \eqref{d19j+3 mod19 1}  follow analogously from Lemma \ref{Lemma Radu 1}  and the chart below.
\end{proof}

\begin{longtable}{ p{.18\textwidth}  p{.40\textwidth}  p{.20\textwidth}  p{.10\textwidth} } 
 \hline
 Congruence  & $(m,M,N,t,(r_\delta))$ and $(r_\delta^\prime)$ & $P(t)$ &
$\floor*{\nu}$ \T\B\\  \hline
\eqref{d2 mod5 1} &
$(125,10,10,97,(3,2,-2,0))$ & \{97,122\} &   22\T\B\\
 & and (30,0,0,0) & &  \T\B\\
 \eqref{d1 mod7 1} &
$(49,14,14,45,(3,1,-1,0))$ & \{45\} &   5\T\B\\
 & and (4,0,0,0) & &  \T\B\\
 & $(49,14,14,17,(3,1,-1,0))$ & \{17,31,38\} &   6\T\B\\
 & and (4,0,0,0) & &  \T\B\\
\eqref{d3 mod49 1} &
$(49,14,14,41,(4,3,-2,0))$ & \{41\} &   12\T\B\\
 & and (9,0,0,0) & &  \T\B\\
\eqref{d3 mod7 1} &
$(343,14,14,90,(39,3,-7,0))$ & \{90,188,237\} &   92\T\B\\
 & and (60,0,0,0) & &  \T\B\\ 
\eqref{d4 mod7 1} &
$(343,14,14,39,(1,4,-2,0))$ & \{39,235,284\} &   76\T\B\\
 & and (77,0,0,0) & &  \T\B\\ 
\eqref{d4 mod11 1} &
$(121,22,22,96,(9,4,-2,0))$ & \{96\} &   31\T\B\\
 & and (11,0,0,0) & &  \T\B\\  
\eqref{d5 mod11 1} &
$(121,22,22,91,(6,5,-2,0))$ & \{91\} &   33\T\B\\
 & and (14,0,0,0) & &  \T\B\\
\eqref{d7 mod11 1} &
$(121,22,22,81,(0,7,-2,0))$ & \{81\} &   34\T\B\\
 & and (19,0,0,0) & &  \T\B\\
\eqref{d6 mod17 1}  &
$(289,34,34,205,(15,6,-2,0))$ & \{205\} &   77\T\B\\
 & and (16,0,0,0) & &  \T\B\\  
 &$(289,34,34,52,(15,6,-2,0))$ & \{52,69,137,171\} &   77\T\B\\
 & and (16,0,0,0) & &  \T\B\\  
 &$(289,34,34,52,(15,6,-2,0))$ & \{188,222,239,273\} &   77\T\B\\
 & and (16,0,0,0) & &  \T\B\\  
\eqref{d19j+3 mod19 1}  &
$(19,38,38,16,(9,3,-1,0))$ & \{16\} &   29\T\B\\   
 & and (1,0,0,0) & &  \T\B\\      
\hline
\end{longtable}

\section*{Acknowledgement}
 The third author was partially supported by Council of Scientific \& Industrial Research (CSIR), Government of India under CSIR-JRF scheme. The author thanks the funding agency.

\end{document}